\documentclass[12pt,a4paper]{amsart}
\usepackage{amssymb}
\usepackage{mathrsfs}
\headsep=1cm \numberwithin{equation}{section}
\newtheorem{thm}{Theorem}[section]
\newtheorem{cor}[thm]{Corollary}
\newtheorem{lem}[thm]{Lemma}
\newtheorem{prop}[thm]{Proposition}
\theoremstyle{definition}

\theoremstyle{question}
\theoremstyle{remark}

\numberwithin{equation}{section}

\newcommand\Ass{\operatorname{Ass}}

\newcommand\Ann{\operatorname{Ann}}
\newcommand\Spec{\operatorname{Spec}}
\newcommand\Rad{\operatorname{Rad}}
\newcommand\Hom{\operatorname{Hom}}
\newcommand\Ext{\operatorname{Ext}}

\newcommand\height{\operatorname{height}}

\newcommand\m{\operatorname{\frak m}}
\newcommand\Soc{\operatorname{Soc}}

\begin{document}\title[Socle finiteness of local cohomology and Gorenstein ideals]
{Socle finiteness of local cohomology modules and Gorenstein ideals}
\author[Ali Akbar Mehrvarz, Kamal Bahmanpour and Reza Naghipour]{Ali Akbar Mehrvarz, Kamal Bahmanpour and Reza Naghipour$^*$}
\address{Department of Mathematics, Islamic Azad University-Tabriz
Branch, Tabriz 5157944533, Iran.}
\address{Department of Mathematics, Islamic Azad University-Ardabil
branch, P.O. Box 5614633167, Ardebil, Iran.}
\email{bahmanpour.k@gmail.com}
\address{Department of Mathematics, University of Tabriz, Tabriz, Iran;
and School of Mathematics, Institute for Research in Fundamental
Sciences (IPM), P.O. Box 19395-5746, Tehran, Iran.}
\email{naghipour@ipm.ir} \email {naghipour@tabrizu.ac.ir}

\thanks{ 
$^*$Corresponding author: e-mail: {\it naghipour@ipm.ir} (Reza Naghipour)}%

\begin{abstract}
The purpose of this paper is to give some equivalent conditions to the
socle and Bass numbers' conjectures which raised by C. Huneke in (Problems on local cohomology,
Free resolutions in commutative algebra and algebraic geometry,
Res. Notes Math. 1992, pp. 93-108). In addition,
some results about certain Gorenstein ideals are included.
\end{abstract}
\maketitle
\section{Introduction}
Let $R$ denote a commutative Noetherian ring
(with identity) and $I$ an ideal of $R$. For an $R$-module $M$, the
$i^{th}$ local cohomology module of $M$ with respect to $I$ is
defined as$$H^i_I(M) = \underset{n\geq1} {\varinjlim}\,\,
\text{Ext}^i_R(R/I^n, M).$$ We refer the reader to \cite{BS} and \cite{Gr1}
 for more details about local cohomology. \\

In \cite{Hu} Huneke conjectured the following:\\

{\bf Conjecture $\rm(i)$:} {\it For any ideal $I$ in a regular local ring
$(R, \frak m)$, the socle of $H^i_I(R)$ is finitely generated for each integer $i\geq 0$.} \\

{\bf Conjecture $\rm(ii)$:} {\it For any ideal $I$ in a regular local ring
$(R, \frak m)$, the Bass numbers
$$\mu^j({\frak p},H^i_I(R))={\rm dim}_{k({\frak p})}\,{\rm
Ext}^j_{R_{\frak p}}(k({\frak p}),H^i_{IR_{\frak p}}(R_{\frak p}))$$
are finite for all integers $i$ and $j$ and all prime ideals ${\frak
p}$ of $R$, where $k({\frak p}):=R_{\frak p}/{\frak p}R_{\frak p}$.
In particular the injective resolution of the local cohomology has
only finitely many copies of the injective hull of $R/{\frak p}$ for
any ${\frak p}$.}\\
There were evidence that the conjectures are true. It is shown by
Huneke and Sharp \cite{HS} and Lyubeznik \cite{Ly1,Ly2} that these two
conjectures hold for a regular ring containing a field. These  problems
are not true in general. For example, if
$I:=(x,y)R\subseteq R:=k[x,y,z,w]/(xz-yw)$, then $\mu^0({\frak
m},H^2_I(R))=\infty$ for ${\frak m}:=(x,y,z,w)$, (cf. \cite{Ha}).

In \cite{AB} it has shown that two conjectures are equivalent.
In Section 2, we give some equivalent conditions to these conjectures. More precisely, we
prove the following result:
\begin{thm}
Let $(R,\m)$ be a  regular local ring. Then the following statements are equivalent:

$\rm(i)$ The Bass numbers $\mu^j(\frak m, H^i_I(R))$ are finite, for each ideal $I$ of $R$ and
for all integers $i, j\geq 0$.

$\rm(ii)$ The $R$-module $\Soc_R(H^j_I(R))$ is finitely generated, for
each
ideal $I$ of $R$ and each integer $j\geq 0$.

$\rm(iii)$ The $R$-module $\Soc_R(H^3_{(x_1,x_2,y)}(R))$ is finitely
generated, whenever $x_1,x_2\in R$ is an $R$-regular sequence and $y\in
R$.

$\rm(iv)$ The $R$-module $\Ext^2_R(R/\m,H^2_{(x_1,x_2,y)}(R))$ is
finitely generated, whenever $x_1,x_2,y\in R$ and ${\rm
grade}(x_1,x_2,y)=2$.
\end{thm}

The proof of Theorem 1.1 is given in 2.6. One of our tools for proving Theorem 1.1 is the following:
\begin{prop}
Let $(R,\m)$ be a regular local ring  and let $x_1,x_2, y$ be elements of $R$
such that $x_1,x_2$ is an $R$-regular sequence. Then the $R$-module $\Soc_R(H^3_{(x_1,x_2,y)}(R))$ is finitely
generated if and only if the $R$-module $\Ext^2_R(R/\m,H^2_{(x_1,x_2,y)}(R))$ is
finitely generated.
\end{prop}

In Section 3, we shall be concerned with the description of those certain ideals $I$
in a Gorenstein local ring $R$, for which the factor ring $R/I$ is also Gorenstein.
These ideals are called Gorenstein ideals (cf. \cite[Definition 2.3]{BDNS}). It is
well known that every ideal $I$ generated by a regular sequence, in a Gorenstein local
ring $R$, is a  Gorenstein ideal (cf. \cite[Ex. 18.1]{Mat}). Especially, in this section, we shall provide
a characterization of Gorenstein principal ideals generated by zero-divisors. A typical
result in this direction is the following:
\begin{thm}
Let $(R, {\frak m})$ be a  Gorenstein local ring and  let $y\in
Z_R(R)$. Then $\Hom_R(R/Ry,R)\cong R/Ry$ if and only if $Ry$ is a Gorenstein ideal.
\end{thm}

Pursuing this point of view further we establish some results about
Gorenstein ideals. In fact, we derive the following consequence of
Theorem 1.3.

\begin{cor}
Let $(R, {\frak m})$ be a local (Noetherian) ring of dimension $d\geq
1$, and $\Ass_RR=\{{\frak p}\}$. Then, if  $\frak p \neq 0$ and principal,
then $(0:_R{\frak p})\cong R/{\frak p}$ . In particular, if $R$ is
Gorenstein, then $\frak p$ is a Gorenstein ideal.
\end{cor}

Throughout this paper, $R$ will always be a commutative Noetherian
ring with non-zero identity. The set
of zerodivisors in $R$ will be denoted by $Z_R(R)$.
For each $R$-module $L$, we denote by
 ${\rm Ass }_RL$ the set of associated prime ideals of
$L$.  Also, for any ideal $\frak{b}$ of $R$, {\it the
radical of} $\frak{b}$, denoted by $\Rad(\frak{b})$, is defined to
be the set $\{x\in R \,: \, x^n \in \frak{b}$ for some $n \in
\mathbb{N}\}$. Finally, if $(R, \frak m)$ is local and $M$ is an $R$-module,
then {\it the  socle of} $M$ is denoted by ${\rm Soc}_R\,(M)$ and defined as $\Hom_R(R/\frak m, M)$.
For any unexplained notation and terminology we refer the reader to \cite{BS}, \cite{BH} and \cite{Mat}.\\

\section{Bass numbers and socle finiteness of local cohomology}

The purpose of this section is to give some equivalent conditions to the
socle and Bass numbers' conjectures which raised by C. Huneke in \cite{Hu}.
The following two lemmas are quite useful in this section.

\begin{lem}
Let $(R,\m)$ be a  Cohen-Macaulay local ring. Then the following two
statements are equivalent:

$\rm(i)$ The $R$-module $\Soc_R(H^j_I(R))$ is finitely generated, for
each ideal $I$ of $R$ and each integer $j\geq 0$.

$\rm(ii)$ The following two conditions are fulfilled:

$\rm (a)$ The $R$-module $\Soc_R(H^2_{(x,y)}(R))$ is finitely generated
for every $x,y\in R$.

$\rm(b)$ The $R$-module $\Soc_R(H^3_{(x_1,x_2,y)}(R))$ is finitely
generated, whenever $x_1,x_2\in R$ is an $R$-regular sequence and $y\in
R$.
\end{lem}
\proof

The assertion follows by the method used in the proof of
\cite[Theorem 6]{He}. \qed\\

\begin{lem}
Let $(R,\m)$ be a regular local  ring. Then the $R$-module
$\Soc_R(H^2_{(x,y)}(R))$ is finitely generated for
every $x,y\in R$.
\end{lem}

\proof We may assume that $H^2_{(x,y)}(R)\neq 0$. Now if ${\rm
heigth}(x,y)=2$, then we have ${\rm grade}(x,y)=2$ and so the
assertion follows from \cite[Theorem 6.2.7]{BS} and \cite[Corollary
3.5]{Kh}. But, if ${\rm heigth}(x,y)=1$, then by \cite[Proposition
2.20]{BN1}, the Bass numbers of the $R$-module $H^1_{(x,y)}(R)$ are
finite. Also, as $R$ is a domain it follows that $H^0_{(x,y)}(R)=0$.
Therefore, in view of \cite[Corollary 3.5]{Kh} the assertion holds.
\qed\\
\begin{prop}
Let $(R,\m)$ be a regular local ring. Then the following two
statements are equivalent:

$\rm(i)$ The $R$-module $\Soc_R(H^j_I(R))$ is finitely generated, for
each ideal $I$ of $R$ and each integer $j\geq 0$.

$\rm(ii)$ The $R$-module $\Soc_R(H^3_{(x_1,x_2,y)}(R))$ is finitely
generated, whenever $x_1,x_2\in R$ is an $R$-regular sequence and $y\in
R$.
\end{prop}
\proof The assertion follows easily from Lemmas 2.1 and 2.2.\qed\\

The following proposition and its corollary are needed in the proof of
the main result of this section.
\begin{prop}
Let $(R,\m)$ be a regular local ring  and let $x_1,x_2, y$ be elements of $R$
such that $x_1,x_2$ is an $R$-regular sequence. Then the $R$-module $\Soc_R(H^3_{(x_1,x_2,y)}(R))$ is finitely
generated if and only if the $R$-module $\Ext^2_R(R/\m,H^2_{(x_1,x_2,y)}(R))$ is
finitely generated.
\end{prop}
\proof Since ${\rm grade}(x_1,x_2,y)\geq 2$, it follows from
\cite[Theorem 6.2.7]{BS} that
$$H^0_{(x_1,x_2,y)}(R)=0=H^1_{(x_1,x_2,y)}(R).$$ Therefore, it follows
from \cite[Corollary 3.5]{Kh}, that the $R$-modules
$$\Ext^i_R(R/\m,H^2_{(x_1,x_2,y)}(R)),$$ are finitely generated for
$i=0,1$. Consequently, in view of \cite[Corollary 3.5]{Kh}, the $R$-module
$\Soc_R(H^3_{(x_1,x_2,y)}(R))$ is finitely generated if and only
if  $\Ext^2_R(R/\m,H^2_{(x_1,x_2,y)}(R))$ is a finitely
generated $R$-module, as required.\qed\\

\begin{cor}
Let $(R,\m)$ be a regular local ring. Then the following two
statements are equivalent:

$\rm(i)$ The $R$-module $\Soc_R(H^j_I(R))$ is finitely generated, for
each ideal $I$ of $R$ and each integer $j\geq 0$.

$\rm(ii)$  The $R$-module $\Ext^2_R(R/\m,H^2_{(x_1,x_2,y)}(R))$ is
finitely generated, whenever $x_1,x_2,y\in R$ such that ${\rm grade}(x_1,x_2,y)=2$.
\end{cor}
\proof It is easy to see that we may assume $x_1\neq0$. Then, as $R$ is domain, $x_1\not \in Z_R(R)$.
Also, since ${\rm grade}(x_1,x_2,y)=2$, it follows that $$Rx_2+(x_1, y)\nsubseteq Z_R(R/Rx_1),$$
and so,  by \cite[Ex. 16.8]{Mat}, there exists $z\in(x_1, y)$ such that $x_2+z\not\in Z_R(R/Rx_1)$.
Now, as $x_1,x_2+z$ is an $R$-regular sequence
and $(x_1,x_2,y)=(x_1,x_2+z,y)$, the assertion follows from the Propositions 2.3 and 2.4.\qed\\

We are now ready to state and prove the main result of this section, which give us some
equivalent conditions to socle and Bass numbers' conjectures (cf. \cite{Hu}).

\begin{thm}
Let $(R,\m)$ be a  regular local ring. Then the following statements are equivalent:

$\rm(i)$ The Bass numbers $\mu^j(\frak m, H^i_I(R))$ are finite, for each ideal $I$ of $R$ and
for all integers $i, j\geq 0$.

$\rm(ii)$ The $R$-module $\Soc_R(H^j_I(R))$ is finitely generated, for
each
ideal $I$ of $R$ and each integer $j\geq 0$.

$\rm(iii)$ The $R$-module $\Soc_R(H^3_{(x_1,x_2,y)}(R))$ is finitely
generated, whenever $x_1,x_2\in R$ is an $R$-regular sequence and $y\in
R$.

$\rm(iv)$ The $R$-module $\Ext^2_R(R/\m,H^2_{(x_1,x_2,y)}(R))$ is
finitely generated, whenever $x_1,x_2,y\in R$ and ${\rm
grade}(x_1,x_2,y)=2$.
\end{thm}
\proof The assertion follows from \cite[Theorem 2.3]{AB}, Proposition 2.4
and Corollary 2.5.\qed\\

\section{Some results on the Gorenstein ideals}

We shall be concerned in this section with the description of those certain ideals $I$
in a Gorenstein local ring $R$, for which the factor ring $R/I$ is also Gorenstein.
These ideals is called Gorenstein ideals (cf. \cite[Definition 2.3]{BDNS}). It is
well known that every ideal $I$ generated by a regular sequence in a local Gorenstein
ring $R$, is a  Gorenstein ideal (cf. \cite[Ex. 18.1]{Mat}). Especially, in this section, we shall provide
a characterization of Gorenstein principal ideals with generated by zero-divisors. Before we state the
the main theorem of this section, we recall a couple lemmas that we will used in the proof of that theorem.

\begin{lem} Let $(R, \frak m)$ be a Cohen-Macaulay local ring.
Then $R$ is Gorenstein if and only if the canonical module $\omega _R$ of $R$
exists and $\omega _R\cong R$.
\end{lem}
\proof See \cite[Theorem 3.3.7]{BH}.\qed\\

\begin{lem} Let $(R, \frak m)$ be a Gorenstein local ring and let $I\subseteq R$ be
an ideal of grade $t$ such that $R/I$ is a Cohen-Macaulay local ring. Let ${\bf x}= x_1, \dots, x_n$
be a system of generators of $I$. Then the canonical module $\omega _{R/I}$ of $R/I$
exists and $\omega _{R/I}\cong H_{n-t}(\bf x)$.
\end{lem}
\proof See \cite[Ex. 3.2.24]{BH}.\qed\\

 We are now ready to state and prove the main theorem of this section, which gives a
 characterization of Gorensteinness of a principal ideal generated by a zerodivisor.\\

\begin{thm}
Let $(R, {\frak m})$ be a  Gorenstein local ring of dimension $n$ and let $y\in
Z_R(R)$. Then the following conditions are equivalent: \\

{\rm(i)}  $\Hom_R(R/Ry,R)\cong R/Ry$.

{\rm(ii)}  $Ry$ is  Gorenstein ideal.
\end{thm}
\proof
In order to show  ${\rm(i)}\Longrightarrow{\rm(ii)}$, by assumption (i), we have
$R/Ry\cong (0:_R Ry)$. Hence there exist the following exact sequences,
$$0\longrightarrow R/Ry\longrightarrow R \longrightarrow Ry\longrightarrow0,$$ and
$$0\longrightarrow Ry\longrightarrow R \longrightarrow R/Ry \longrightarrow0.$$
As $R$ is Cohen-Macaulay, it follows that $H^i_{\frak m}(R)=0$ for all $0\leq i\leq n-1$.
Combining this with the long exact sequences of local cohomology modules induced by the above
short exact sequences, we have the following
\begin{center}
$H^0_{\frak m}(Ry)=H^0_{\frak m}(R/Ry)=0,$\vspace*{0.2cm}

$H^i_{\frak m}(R/Ry)\cong H^{i+1}_{\frak m}(Ry),$  for all $i\leq n-2,$\vspace*{0.2cm}

$H^j_{\frak m}(Ry)\cong H^{j+1}_{\frak m}(R/Ry), $ for all $j\leq n-2$.
\end{center}

It follows that $H^i_{\frak m}(R/Ry)=0$ for all $i\leq n-1.$  Hence $R/Ry$ is Cohen-Macaulay.
As $R$ is Gorenstein it follows from Lemma 3.2 that $\omega _{R/Ry}\cong H_0(y)\cong R/Ry$.
Now the result follows from Lemma 3.1.

The implication ${\rm(ii)}\Longrightarrow{\rm(i)}$  follows from Lemma 3.2.  \qed\\
\begin{cor}
Let $(R, {\frak m})$ be a local (Noetherian) ring of dimension $d\geq
1$, and $\Ass_RR=\{{\frak p}\}$. Then, if  $\frak p \neq 0$ and principal,
then $(0:_R{\frak p})\cong R/{\frak p}$ . In particular, if $R$ is
Gorenstein, then $\frak p$ is a Gorenstein ideal.
\end{cor}
\proof It follows from $\Ass_RR=\{{\frak p}\}$ that ${\frak p}={\rm nil}(R)$.
As ${\frak p}\neq
0$, it follows that there is a positive integer $n$ such that
${\frak p}^n\neq 0$ and ${\frak p}^{n+1}=0$. Let $\frak p=Rx$
for some $0\neq x\in R$. Then $$(0:_R{\frak
p})\subseteq Z_R(R)={\frak p}=Rx,$$
and so there is an ideal $\frak a$ of $R$ such that $(0:_R{\frak p})=\frak a x$.
Now, as $r\in \frak a$ if and only if $rxx=0$ if and only if $r\in (0:_R{\frak p}^2)$,
it follows that $(0:_R{\frak p})=(0:_R{\frak p}^2)x$. Similarly, it yields that
$(0:_R{\frak p}^2)=(0:_R{\frak p}^3)x$, and so we have
$(0:_R{\frak p})=(0:_R{\frak p}^3)x^2$. Hence, by induction, we obtain
$$(0:_R{\frak p})=(0:_R{\frak p}^n)x^{n-1}.$$
Furthermore, it is easy to see that $(0:_R{\frak p}^n)={\frak p}$, and so
$(0:_R{\frak p})=Rx^n$. Therefore,
$$(0:_R{\frak p})=Rx^n\cong R/(0:_Rx^n)= R/(0:_R \frak p^n)=R/\frak p.$$

Now, the second assertion follows from Theorem 3.3. \qed\\

\begin{thm}
Let $(R, {\frak m})$ be a  Gorenstein local ring. Let ${\frak p}\in \Spec R$
and let $y_1, \dots, y_t\in R$ be an $R$-regular sequence. Let $x\in
R$ be such that $\frak p=(y_1,\dots, y_t, x)$ and
$\Rad(y_1, \dots, y_t)=\frak p$. Then $\frak p$ is a Gorenstein ideal.
\end{thm}

\proof If ${\frak p}=(y_1,\dots, y_t)$, then the result follows from
\cite[Ex. 18.1]{Mat}. So we can assume that ${\frak p}\neq (y_1,\dots, y_t)$.
 Then $x\not\in (y_1,\dots, y_t)$. In order to simplify notation, we will
use $\bar{R}$ to denote the  Gorenstein local ring
$R/(y_1,\dots, y_t)$ (see \cite[Ex. 18.1]{Mat}) and we will write
$\bar{x}$ for the element $x+ (y_1,\dots, y_t)$ of
$\bar{R}$. Then, as $\Rad(y_1, \dots, y_t)=\frak p,$ it follows
that $Ass_{\bar{R}}\bar{R}=\{\frak p/(y_1,\dots, y_t)\}$. Now,
the result follows from Corollary 3.4.  \qed\\

\begin{thm} Let $(R, {\frak m})$ be a Gorenstein local ring. Let
$x_1,\dots,x_n$ be an $R$-regular sequence and $y\in {\frak m}$ such that $\height(x_1,\dots,x_n,y)=n$.
Then the following conditions are equivalent:\\

{\rm (i)} $\Ext^n_R(R/(x_1,\dots,x_n, y), R)\cong R/(x_1,\dots,x_n,y)$.

{\rm (ii)} $(x_1,\dots,x_n, y)$ is a Gorenstein ideal.
\end{thm}
\proof In order to simplify notation, we will
use $\bar{R}$ to denote the  Gorenstein local ring
$R/(x_1,\dots, x_n)$ (see \cite[Ex. 18.1]{Mat}) and we will write
$\bar{y}$ for the element $y+ (x_1,\dots, x_n)$ of
$\bar{R}$. Then, an application of Rees' Theorem (cf. \cite[Lemma 1.2.4]{BH}) gives
the isomorphism
$$\Ext^n_R(R/(x_1,\dots,x_n, y), R)\cong \Hom_{\bar{R}}(\bar{R}/\bar{y}\bar{R}, \bar{R}).$$
Moreover, as  $\height(x_1,\dots,x_n,y)=n$, it follows that $\bar{y}\in Z_{\bar{R}}(\bar{R})$.
Now, the result follows from the Theorem 3.3 and \cite[Ex. 18.1]{Mat}. \qed\\

The following lemma, which is a consequence of the local duality theorem for Gorenstein local rings,
and is needed in the proof of the Theorem 3.8, is almost certainly know, but we could not find
a reference for it, so it is explicit stated and proved here.

\begin{lem} Let $(R, \frak m)$ be a Gorenstein local ring and let $I$ be
an ideal of $R$ such that $\dim R/I=\dim R$. Then the following conditions are equivalent:

{\rm(i)} $\Hom_R(R/I,R)\cong R/I$ and $\Ext^i_R(R/I,R)=0$ for each $i\geq 1$.

{\rm(ii)} $I$ is  Gorenstein ideal.
\end{lem}
\proof  Since both conditions (i) and (ii) behave well under completion, without loss of
generality, in view of \cite[Theorem 18.3]{Mat}, we may assume that $(R, \frak m)$ is complete.

To prove the implication ${\rm(i)}\Longrightarrow{\rm(ii)}$, let $k:=R/\frak m$,  $E:=E(k)$ and
$D:=\Hom_R(\cdot, E)$. Then, by \cite[Theorem 11.2.5]{BS}, we have
$$H^{d-i}_I(R/I)\cong D(\Ext^i_R(R/I, R)),$$
for all integers $i$, where $d:=\dim R$. Thus by assumption (i) and \cite[Theorem 10.1.15]{BS},
we have $H^{d}_I(R/I)\cong E_{R/I}(k)$ and $H^{i}_I(R/I)=0$ for $i\neq d$. Therefore it follows from
\cite[Corollary 9.5.13]{EJ} that $I$ is Gorenstein ideal.

In order to prove ${\rm(ii)}\Longrightarrow{\rm(i)}$, as $R$ is complete, using local duality theorem
for Gorenstein local ring \cite[Theorem 11.2.5]{BS}, we obtain
$$D(H^{d-i}_I(R/I))\cong \Ext^i_R(R/I, R),$$ for all integers $i$. Hence, again using \cite[Corollary 9.5.13]{EJ},
it is easy to see that $\Hom_R(R/I,R)\cong R/I$ and $\Ext^i_R(R/I,R)=0$, for all $i\geq 1$. \qed\\

\begin{thm}
Let $(R, \frak m)$ be a  Gorenstein local ring and $I$  a
non-nilpotent ideal of $R$ such that $\dim R/I=\dim R$. Then $I^n$ is not Gorenstein ideal, for all large $n$.
\end{thm}
\proof There exists a positive integer $k$ such that
$$(0:_RI^k)=(0:_RI^{k+1})=\cdots.$$
Now, for $n\geq k+1$,  we have $$I^k\subseteq
\Ann_R(0:_RI^k)=\Ann_R(0:_RI^n)=\Ann_R(\Hom_R(R/I^n,R)),$$ and so
$\Hom_R(R/I^n,R)\not\cong R/I^n.$
Therefore, it follows from Lemma 3.7 that the ideal $I^n$ is not Gorenstein, for all $n\geq
k+1$.\qed

\begin{cor}
Let $(R, \frak m)$ be a  Gorenstein local ring and $y$ a
non-nilpotent element of $R$ such that $y\in Z_R(R)$. Then, the ideal $Ry^n$ is
not Gorenstein, for all large $n$.
\end{cor}
\proof The assertion follows from Theorem 3.8. (Note that $\dim R/Ry=\dim R.$) \qed\\

{\bf Remark.} Let $R$ be a regular local ring, and let $x_1, x_2, y\in R$ be such that $x_1, x_2$ is an $R$-regular sequence.
Then, in view of Theorem 2.6, it is shown that the $R$-module ${\rm Soc}_R(H^3_{(x_1,x_2,y)}(R))$ is finitely generated
if and only if the $R$-module $\Soc_R(H^j_I(R))$ is finitely generated, for
each ideal $I$ of $R$ and each integer $j\geq 0$. Unfortunately, however we have not been able
to prove that the $R$-module $\Soc_R(H^j_I(R))$ is finitely generated. On the other hand, it is known that
$H^3_{(x_1,x_2,y)}(R)\cong \underset{n\geq1} {\varinjlim}\,\,R/(x^n_1, x^n_2, y^n),$
and so we end this paper with the following question.\\

{\bf Question.} Is the set
$\Omega:=\{\dim _{R/\frak m}{\rm Soc}_R(R/(x^n_1,x^n_2,y^n))| \,\, n\in\Bbb{N}\}$ bounded from
above?\\


\begin{center}
{\bf Acknowledgments}
\end{center}

 The authors are deeply grateful to the
referee for his/her careful reading and many helpful suggestions on
the paper. Also, the authors would like to thank the Islamic Azad University-Tabriz Branch for its financial support.

\end{document}